\documentclass[openany, a4paper, 11pt]{article}

\usepackage{amsfonts}

 \usepackage{mathrsfs,amsfonts,amsmath}
 \usepackage{color}

 \makeatletter
 \@addtoreset{equation}{section}
 \makeatother
 \newtheorem{theorem}{Theorem}[section]
 \newtheorem{definition}{Definition}[section]
 \newtheorem{hypothesis}{Hypothesis}[section]
 \newtheorem{lemma}{Lemma}[section]
 \newtheorem{proposition}{Proposition}[section]
 \newtheorem{corollary}{Corollary}[section]
 \newtheorem{remark}{Remark}[section]
 \newtheorem{example}{Example}[section]

 \def\beqlb{\begin{eqnarray}}\def\eeqlb{\end{eqnarray}}
 \def\beqnn{\begin{eqnarray*}}\def\eeqnn{\end{eqnarray*}}

 \def\ar{\!\!&}

 \def\mbb{\mathbb}

 \def\qed{\hfill$\Box$\medskip}
 \def\dfR{{\mbb R}}

\newcommand{\bcen}{\begin{center}}
\newcommand{\ecen}{\end{center}}
\newcommand{\bgeqn}{\begin{equation}}
\newcommand{\edeqn}{\end{equation}}

\def\dz{\delta}

\def\ez{\epsilon}

\def\D{{\cal D}}
\def\E{{\bf E}}
\def\G{{\cal G}}
\def\L{{\cal L}}

\def\L{{\cal L}}

\def\S{{\cal S}}
\def\tlb{{\tilde{B}}}

\def\rar{\rightarrow}

\def\ra{\rangle}
\def\la{\langle}

\begin{document}

\centerline{\Large\textbf{ Fleming-Viot Processes in an Environment
}\footnote{Supported by NSFC (No.10721091 )}}

\bigskip

\centerline{ Hui He\footnote{ \textit{E-mail address:} {
hehui@bnu.edu.cn }} }

\medskip

\centerline{Laboratory of Mathematics and Complex Systems, }

\centerline{ School of Mathematical Sciences, Beijing Normal
University,}

\smallskip

\centerline{ Beijing 100875, People's Republic of China}

\bigskip

{\narrower{\narrower{\narrower

\begin{abstract}
We consider a new type of lookdown processes where spatial motion of
each individual is influenced by an individual noise and a common
noise, which could be regarded as an environment. Then a class of
probability measure-valued processes on real line $\mbb{R}$ are
constructed.
 The sample path properties are investigated: the values of
this new type process are either purely atomic measures or
absolutely continuous measures according to the existence of
individual noise. When the process is absolutely continuous with
respect to Lebesgue measure, we derive a new  stochastic partial
differential equation for the density process. At last we show that
such processes also arise from normalizing a class of measure-valued
branching diffusions in a Brownian medium as the classical result
that Dawson-Watanabe superprocesses, conditioned to have total mass
one, are Fleming-Viot superprocesses.
\end{abstract}

\smallskip


\smallskip


\smallskip

\noindent\textbf{Abbreviated Title:} Fleming-Viot processes

\par}\par}\par}

\bigskip\bigskip

\section{Introduction}
In this work, we construct and study a new class of probability
measure-valued Markov processes on the real line $\mbb R$. Our model
arises from a modified stepwise mutation model (see Section 1.1.10
of \cite{[E00]} for classical stepwise mutation model): the mutation
process of each individual in the model is influenced by an
independent noise and a common noise. More precisely, suppose that
$\{ W(t,x): x\in \mathbb{R}, t\geq0\}$ is space-time white noise
based on Lebesgue measure, the common noise, and
$\{B_i(t):t\geq0,i=1,2,\cdots\}$ is a family of independent standard
Brownian motions, the individual noises, which are independent of
$\{W(t,x): x\in\mathbb{R}\}$. The mutation of  an individual in the
stepwise mutation system with label $i$ is defined by the stochastic
equations
 \bgeqn
  \label{1.6}
dx_i(t)=\ez dB_i(t)+\int_\mathbb{R}h(y-x_i(t))W(dt,dy),\textrm{\ \ \
}t\geq0,~~i=1,2,\cdots,
 \edeqn
where $W(dt,dy)$ denotes the time-space stochastic integral relative
to $\{W_t(B)\}$ and $\ez\geq0.$ Suppose that  $h\in C^2(\mathbb{R})$
is square-integrable. Let $\rho_{\ez}=\ez^2+\rho(0)$ and
 \bgeqn \label{f1.1}
\rho(x)=\int_{\mathbb{R}}h(y-x)h(y)dy,
 \edeqn
 for $x\in \mathbb{R}.$   For each integer $m\geq1$,
$\{(x_1(t),\cdots,x_m(t)):t\geq0\}$ is an $m$-dimensional diffusion
process which is generated by the differential operator
\bgeqn\label{Gdiff}
 G^m:=\frac{1}{2}\sum_{i=1}^ma(x_i)\frac{\partial^2}{\partial
x_i^2}+ \frac{1}{2}\sum_{i,j=1,i\neq
j}^m\rho(x_i-x_j)\frac{\partial^2}{\partial x_i\partial x_j}.
 \edeqn
In particular, $\{x_i(t):t\geq0\}$ is a one-dimensional diffusion
process with generator $G:=(\rho_{\ez}/2)\Delta$. Because of the
exchangeability, a diffusion process generated by $G^m$ can be
regarded as an interacting particle system or a measure-valued
process. Heuristically, $\rho_{\ez}$ represents the speed of the
particles and $\rho(\cdot)$ describes the interaction between them.

Our interest comes from recent studies on connections between
superprocesses and stochastic flows; see \cite{[DLW01]},
\cite{[DLZ04]},  \cite{[SA01]} and \cite{[W98]}. In those works,
particles undergo random branching and their spatial motions are
affected by the presence of stochastic flows. Some new classes of
measure-valued processes were constructed from the empirical measure
of the particles. Those measure-valued processes are quite different
with the classical Dawson-Watanabe processes. There are at least two
different ways to look at those processes. One is as a superprocess
in random environment and the other as an extension of models of the
motion of the mass by stochastic flows; see \cite{[LR05]},
\cite{[MX01]}. In this work we remove the branching structure of
particle systems in \cite{[W98]} but add a sampling mechanism. That
is whenever a particle's exponential `sampling clock' rings, it
jumps to a position chosen at random from the current empirical
distribution of the whole population. Its mutation then continues
from its new position.

This work is simulated by classical connections between
Dawson-Watanabe processes and Fleming-Viot processes investigated in
\cite{[EM91]} and \cite{[P92]}. It has been shown that  Fleming-Viot
superprocesses is the Dawson-Watanabe prcesses, conditioned to have
total mass one. So we want to ask what can we obtain if the
measure-valued processes constructed in \cite{[DLW01]},
\cite{[DLZ04]} \cite{[SA01]} and \cite{[W98]} are conditioned to
have total mass one? The particle picture described in \cite{[P92]}
suggests that the branching structure of such conditioned
measure-valued branching processes may be changed to sampling
mechanism. Thus measure-valued branching processes constructed in
\cite{[W98]}, conditioned to have total mass one, may have generator
as:
 \bgeqn
 \label{f1.2}
 \mathcal{L}F(\mu):=\mathcal{A}F(\mu)+\mathcal{B}F(\mu),
 \edeqn
where
 \begin{eqnarray} \label{f1.3}
\mathcal{A}F(\mu)\ar:=\ar\frac{1}{2}\int_{\mathbb{R}}\rho_{\ez}\frac{d^2}{dx^2}
                    \frac{\dz F(\mu)}{\dz\mu(x)}\mu(dx)\cr\ar\ar
+\frac{1}{2}\int_{\mathbb{R}^2}\rho(x-y)\frac{d^2}{dxdy}
                  \frac{\dz^2F(\mu)}{\dz\mu(x)\dz\mu(y)}\mu(dx)\mu(dy),
\end{eqnarray}
 \bgeqn \label{f1.4}
\mathcal{B}F(\mu):=\frac{\gamma}{2}\int_\mathbb{R}\int_{\mbb
R}\frac{\dz^2F(\mu)}{\dz\mu(x)\dz\mu(y)}
 \left(\mu(dx)\dz_{x}(dy)-\mu(dx)\mu(dy)\right),
 \edeqn
for some bounded continuous functions $F(\mu)$ on $P(\mathbb{R})$.
The variational derivative is defined by
\begin{equation}
\label{f1.5} \frac{\dz
F(\mu)}{\dz\mu(x)}=\lim_{r\rightarrow{0+}}\frac{1}{r}[F(\mu+r\dz_x)-F(\mu)],\textrm{\
\ \ } x\in \mathbb{R},
\end{equation}
if the limit exists and  $\dz^2F(\mu) / \dz\mu(x)\dz\mu(y)$ is
defined in the same way with $F$ replaced by $(\dz F/\dz\mu(y))$ on
the right hand side. If we replace $\cal B$ in (\ref{f1.4}) by
$$
\frac{\gamma}{2}\int_\mathbb{R}
\frac{\dz^2F(\mu)}{\dz\mu(x)^2}\mu(dx),
$$
then $\cal L$ is  the generator of the measure-valued process
constructed in \cite{[W98]}, where $\cal L$ acted on some bounded
continuous functions on $M(\mbb R)$, space of finite measures on
$\mbb R$; see (1.8) of \cite{[W98]}. If the second term in $\cal A$
vanishes, then $\cal L$ is just the generator of an usual
Fleming-Viot process.
\smallskip

The main work in this paper is to solve the martingale problem and
analyze the sample path properties of the solution.  For $f\in
B(\mbb R^m)$, define $F_{m,f}(\mu)=\langle f,\mu^m\rangle$.
 For
$\mu\in P(\mathbb{R})$, we say a $P(\mathbb{R})$-valued continuous
process $\{Z(t):t\geq0\}$ is a solution of the $(\mathcal
{L},\mu)$-\textsl{martingale problem} if $Z(0)=\mu$ and
 \beqlb\label{martfor}
  F(Z(t))-F(Z(0))-\int_0^t\mathcal {L}F(Z(s))ds,\textrm{\ \ \ } t\geq0,
 \eeqlb
is a martingale for each $F\in\mathcal {D}(\mathcal
{L}):=\bigcup_{m\geq1}\{F_{m,f}(\mu),  f\in C^2(\mathbb{R}^m)\}.$
 A
simple calculation yields
 \beqlb\label{GeneratorL}
\L F_{m,f}(\mu)=\la \mu^m, G^mf\ra+\sum_{1\leq i<j\leq
m}\gamma\left( \la \mu^{m-1},\Psi_{ij}f\ra-\la\mu^m, f\ra\right),
 \eeqlb
where $\Psi_{ij}$ denotes the operator from $B(\mathbb{R}^{m})$ to
$B(\mathbb{R}^{m-1})$ defined by
 \bgeqn
 \label{replace}
 \Psi_{ij}f(x_1,\cdots,x_{m-1})=f(x_1,\cdots,x_{m-1},\cdots,x_{m-1},\cdots,x_{m-2}),
 \edeqn
where $x_{m-1}$ is in the places of the $i$th and the $j$th
variables of $f$ on the right hand side. We shall show that the
$({\cal L}, \mu)$-martingale problem is well-posed and call the
solution as Fleming-Viot process in an environment (FVE for short).
We will use look-down construction suggested by \cite{[DK96]} with
some modifications to show the existence of the solution. This
look-down construction will help us on analyzing the sample path
properties. The uniqueness of the $({\cal L},\mu)$-martingale
problem will be proved by classical duality argument. Since the
spatial motions of individuals in the look-down system are not
independent with each other, when solving the martingale problem, we
need some technical lemmas which will be given in the Appendix.

 Our other
main results include:
\begin{enumerate}

\item State classification: when $\ez>0$,
FVE
 is absolutely
continuous respect to $dx$  and we also deduce a new SPDE for the
density process; when $\ez=0$ its values are purely atomic;

\item When conditioned to have total mass one, a measure-valued branching process
in a Brownian medium constructed in \cite{[W98]} is an FVE.
 \end{enumerate}

The remaining of this paper is organized as follows. In Section 2,
we solve the $({\cal L}, \mu)$-martingale problem. The state
classification of the process will be investigated in Section 3. In
the last section, Section 4, we derive the connection between FVE
and the process constructed in \cite{[W98]}. Two technical lemmas
will be given in the Appendix.
\begin{remark}  By Theorem 8.2.5 of
\cite{[EK86]}, the closure of $\{(f,G^mf):f\in
C_c^{\infty}(\mathbb{R}^m)\}$ denoted by $\bar{G}^m$ is
single-valued and generates a Feller semigroup $(T_t^m)_{t\geq0}$ on
$\hat{C}(\mathbb{R}^m)$. Note that this semigroup is given by a
transition probability function and can therefore be extended to all
of $B(\mathbb{R}^m)$.
\end{remark}
Notation: For reader's convenience, we introduce here our main
notation. Let $\hat{\mathbb{R}}$ denote the  one-point
compactification of $\mathbb{R}$.    Given a topological space $E$,
let $M(E)$ ($P(E)$) denote space of finite measures (probability
measures) on $E$. Let $B(E)$ denote the set of bounded measurable
functions on $E$ and
  let $C(E)$ denote its subset comprising of bounded continuous
  functions. Let $\hat{C}(\mathbb{R}^n)$ be the space of continuous
   functions  on $\mbb R^n$ which vanish at infinity and
  let $C_c^{\infty}(\mathbb{R}^n)$ be functions with compact support
 and bounded continuous derivatives of any order.
  Let $C^2(\mathbb{R}^n)$ denote the set of functions in
$C(\mathbb{R}^n)$ which is twice continuously differential functions
with bounded derivatives up to the second order. Let
$\hat{C}^2(\mathbb{R}^n)$ be the subset of $C^2(\mathbb{R}^n)$ of
functions that together with their
derivatives up to the second order vanish  at infinity. \\
\noindent Let
$$
C_{\partial}^2(\mbb{R}^n)=\{f+c: c\in\mbb{R} \textrm{ and } f\in
\hat{C}^2({\mbb{R}^n})\}
$$
 We denote by $C_E[0,\infty)$ the space of continuous paths taking values in
 $E$.
 Let $D_E[0,\infty)$ denote the Skorokhod space of c\`{a}dl\`{a}g paths taking values
 in $E$. For $f\in C(\mathbb{R})$ and $\mu\in M(\mathbb{R})$ we shall write
$\langle \mu, f\rangle$ for $\int fd\mu$.

\section{Construction}\label{SECCON}
\subsection{Uniqueness}\label{SECUNIQUE}
In this subsection, we define a dual process to show the uniqueness
of the $({\cal L},\mu)$-martingale problem. Let $\{M_t:t\geq0\}$ be
a nonnegative integer-valued c\`{a}dl\`{a}g Markov process. For
$i\geq j$, the transition intensities $q_{i,i-1}=\gamma i(i-1)/2$
and $q_{ij}=0$ for all other pairs $i,j$. Let $\tau_0=0$  and let
$\{\tau_k:1\leq k\leq M_0-1\}$ be the sequence of jump times of
$\{M_t:t\geq0\}$. That is $\tau_1=\inf\{t\geq0:M_t\neq
M_0\},\cdots,\tau_k=\inf\{t>\tau_{k-1}:M_t\neq M_{\tau_{k-1}}\} .$\\
Let $\{\Gamma_k:1\leq k\leq M_0-1\}$ be a sequence of random
operators which are conditionally independent given $\{M_t :t\geq
0\}$ and satisfy
$$\textbf{P}\{\Gamma_k=\Psi_{ij}|M(\tau_k-)=l,M(\tau_k)=l-1\}
=\left(\begin{array}{c}l\cr 2\end{array}\right)^{-1} ,\textrm{\ \ \
}1\leq i< j\leq l,$$
 where $\Psi_{ij}$
 are defined by (\ref{replace}). Let $\textbf{B}$ denote the topological
union of $\{B(\mathbb{R}^m):m = 1,2,\cdots\}$ endowed with pointwise
convergence on each $B(\mathbb{R}^m)$. Then
 \bgeqn
 \label{functionf}F_t={T}_{t-\tau_k}^{M_{\tau_k}}\Gamma_k
 {T}_{\tau_k-\tau_{k-1}}^{M_{\tau_{k-1}}}\Gamma_{k-1}\cdots
 {T}_{\tau_2-\tau_1}^{M_{\tau_1}}\Gamma_1{T}_{\tau_1}^{M_0}F_0, \textrm{\
 \ \ }\tau_k\leq t<\tau_{k+1},~~0\leq k\leq M_0-1,
 \edeqn
defines a Markov process $\{F_t:t\geq0\}$ taking values from ${\bf
B}$. Clearly, $\{(M_t,F_t):t\geq0\}$ is also a Markov process. Let
$\textbf{E}_{m,f}$ denote the expectation given $M_0=m$ and
$F_0=f\in B(\mathbb{R}^m)$.
\begin{theorem}
\label{ThDual}
  Suppose that  $\{Z(t):t\geq 0\}$ is a solution
of the $({\mathcal{L}},\mu)$-martingale problem and assume that
$\{Z(t):t\geq0\}$ and $\{(M_t,F_t):t\geq0\}$ are defined on the same
probability space and independent of each other, then
 \bgeqn
 \label{Dual}
 {\bf{E}}\left\langle Z(t)^m, f\right\rangle
 ={\bf{E}}_{m,f}\big{[}\left\langle
  \mu^{M_t},F_t\right\rangle
 \big{]}
 \edeqn
for any $t\geq0$, $f\in C(\mathbb{R}^m)$ and integer $m\geq1$.
\end{theorem}
\textbf{Proof}. In this proof we set
$F_{\mu}(m,f)=F_{m,f}(\mu)=\langle \mu^m ,f\rangle$.
 It suffices to prove (\ref{Dual}) for $f\in C^2(\mbb R^m)$.
  By the definition of $F_t$ and elementary properties of $M_t$,
we know that $\{(M_t,F_t): t\geq0\}$ has weak generator $\mathcal
{L}^\#$ given by
 \begin{eqnarray}
 \label{ThDuala}
 \mathcal {L}^\#F_{\mu}(m,f)=
 F_{\mu}(m,{G}^mf)+
   \sum_{1\leq i<j\leq
   m}\gamma\left(F_{\mu}(m-1,\Psi_{ij}f)-F_{\mu}(m,f)\right)
 \end{eqnarray}
with $f\in C^2(\mbb{R}^m)$. In view of (\ref{GeneratorL}) we have
 \bgeqn
 \label{ThDualb}
 \mathcal {L}^\#F_{\mu}(m,f)={\mathcal
 {L}}F_{m,f}(\mu).
 \edeqn
 Thus if we can show that for $F_0\in
C^2(\mbb R^m)$, $F_t\in C^2(\mbb R^m)$ for all $t\geq0$, then dual
relationship (\ref{Dual}) follows from Corollary 4.4.13 of
\cite{[EK86]}. To this end, it suffices to show that $T_t^m C^2(\mbb
R^m)\subset C^2(\mbb R^m).$ When $\ez>0$, $G^m$ is uniform elliptic.
The desired result follows from Theorem 0.5 on page 227 of
\cite{[Dy65]}. When $\ez=0$, Lemma \ref{lehomo} yields the desired
conclusion. We are done. \qed

\subsection{Look Down Processes}\label{SECLD}
Suppose that $x_t=(x_1(t),\cdots,x_m(t))$ is a Markov process in
$\mathbb{R}^m$ generated by ${G}^m$. By Lemma 2.3.2 of \cite{[D93]}
we know that $(x_1(t),\cdots,x_m(t))$ is an exchangeable Feller
process. Let $P_t^{(m)}$ denote its transition semigroup. Then
$\{P_t^{(m)}, m\geq1\}$ is a consistent family of Feller semigroups
on $C(\dfR^m)$, i.e., for all $k\leq m$, any $k$-component of
$G^m$-diffusion evolve as a $G^k$-diffusion.

\bigskip

Let $\{B_{ijk},\,1\leq i<j,\,1\leq k<\infty\}$ and $\{B_{i0},\,
i\geq1\}$ be independent Brownian motions, independent of $W$. Let
$\{N_{ij},\,1\leq i<j\}$ be independent, unit rate Poisson
processes, independent of $\{B_{ijk}\},\, W$ and let $\tau_{ijk}$
denote the $k$th jump time of $N_{ij}$. Let $\{X_i(0),\, i\geq1\}$
be an exchangeable sequence of random variables, independent of
$\{U_{ijk}\}$, $\{U_{i0}\}$, $W$ and $\{N_{ij}\}$. Define
$\gamma_{ijk}=\min\{\tau_{i'jk'},\, i'<j:
\tau_{i'jk'}>\tau_{ijk}\}$; that is, $\gamma_{ijk}$ is the first
jump time of $N_j\equiv \sum_{i<j}N_{ij}$ after $\tau_{ijk}$, and
define $\gamma_{j0}=\min\{\tau_{ij1}: i<j\}$. Finally, for $0\leq
 t<\gamma_{j0}$ define
 \beqlb
 \label{lookdown1}
 X_j(t)=X_j(0)
 +\ez B_{j0}(t)+\int_0^t\int_{\dfR}h(y-X_{j}(s))W(dyds)
 \eeqlb
and for $\tau_{ijk}\leq
 t<\gamma_{ijk}$,
 \beqlb
 \label{lookdown2}
 X_j(t)=X_{i}(\tau_{ijk})+\ez (B_{ijk}(t)-B_{ijk}({\tau_{ijk}}))
 +\int_{\tau_{ijk}}^t\int_{\dfR}h(y-X_{j}(s))W(dyds).
 \eeqlb
Since $G^m$-diffusion is an exchangeable consistent family of Feller
diffusions, between the jump times of the Poisson processes, the
$X_j$ behave as a $G^1$-diffusion and any $n$-component of the
particle systems evolve as a $G^n$-diffusion. At the jump times of
$N_{ij}$, $X_j$ ``looks down'' at $X_i$, assumes the value of $X_i$
at the jump time, and then evolves as a $G^1$-diffusion and also any
$n$-component of the particle systems evolve as a $G^n$-diffusion.
Then $X=(X_1,X_2,\cdots)$ is a Markov process with generator given
by
 \beqlb
 \label{generatorX}
{\mbb A}f(x_1,\cdots,x_m)&=& G^mf(x_1,\cdots,x_m)\cr &&+\sum_{1\leq
i<j\leq
m}\left(f(\theta_{ij}(x_1,\cdots,x_m))-f(x_1,\cdots,x_m)\right),
 \eeqlb
where $f\in C^2(\dfR^m)$ and $\theta_{ij}(x_1,\cdots,x_m)$ denote
the element of $\dfR^m$ obtained by replacing $x_j$ by $x_i$ in
$(x_1,\cdots,x_m)$.
\\ \smallskip

As in \cite{[DK96]}, we want to compare the $\dfR^{\infty}$-valued
process $X$ to a sequence of modified Moran-type models. Let $S_m$
denote the collection of permutations of $(1,\cdots,m)$ which we
write as ordered $m$-tuples $s=(s_1,\cdots,s_m)$. Let
$\pi_{ij}:S_m\rar S_m$ denote the mapping such that $\pi_{ij}s$ is
obtained from $s$ by interchanging $s_i$ and $s_j$ and let
$\{M_{ijk}:1\leq i\neq j\leq m,k\geq1\}$ be independent random
mappings $M_{ijk}: S_m\rar S_m$ such that
$P\{M_{ijk}s=s\}=P\{M_{ijk}s=\pi_{ij}s\}=\frac{1}{2}$. In following
we define an $S_m$-valued process $\Sigma^m$ and counting processes
$\{\tilde{N}_{ij}, 1\leq i\neq j\leq m\}$ recursively. Let
$\Sigma^m(0)$ be uniformly distributed on $S_m$ and independent of
all other processes. Let
 \beqlb
 \label{tildeN}
 \tilde{N}_{ij}(t)=\sum_{1\leq k< l\leq m}\int_0^t {\bf1}
 _{\{\Sigma_i^m(r-)=k,\,\Sigma_j^m(r-)=l\}}dN_{kl}(r)
 \eeqlb
and let $\Sigma^m$ be constant except for discontinuities determined
by
$\Sigma^m(\tilde{\tau}_{ijk})=M_{ijk}\Sigma^m(\tilde{\tau}_{ijk}-),$
where $\tilde{\tau}_{ijk}$ is the $k$-th jump time of
$\tilde{N}_{ij}$, or more precisely, interpreting $\Sigma^m$ as a
$\mbb Z^m$-valued process,
 \beqlb
 \label{morans}
 {\Sigma}^{m}(t)=\sum_{1\leq i<j\leq m}\int_0^t
 \left(M_{ij(\tilde{N}_{ij}(r-)+1)}\Sigma^m(r-)\right)
 d\tilde{N}_{ij}(r).
 \eeqlb
Next, define $\{\hat{N}_{ij},1\leq i\leq m<j\}$ by
 \beqlb
 \label{hatN}
 \hat{N}_{ij}(t)=\sum_{k=1}^m\int_0^t{\bf
 1}_{\{\Sigma_i^m(r-)=k\}}dN_{kj}(r) \eeqlb and let
 $\hat{\tau}_{ijk}$ denote the $k$-th jump time of $\hat{N}_{ij}$.
 Note that for $j>m$,
 \beqlb\label{moranN}
 N_j=\sum_{1\leq i<j}N_{ij}=\sum_{1\leq i\leq m}\hat{N}_{ij}+\sum
 _{m<i\leq j}N_{ij}.
 \eeqlb
By Lemma 2.1 of \cite{[DK96]},  $\{\tilde{N}_{ij}\}$ and
$\{\hat{N}_{ij}\}$ are Poisson processes with intensities
$\frac{1}{2}$ and 1, respectively. And for each $t\geq0$,
$\Sigma^m(t)$ is independent of
$\G_t=\sigma(\tilde{N}_{ij}(s),\hat{N}_{kl}(s):s\leq t,1\leq i\neq
j\leq m, 1\leq k\leq m<l)$. Define
$$Y_j^m(t)=X_{\Sigma_j^m(t)}(t),\quad j=1,\cdots,m.$$

\begin{lemma}\label{lemmaY}$
Y^m=(Y_1^m,\cdots, Y_m^m)$ is a Markov process with generator
given by
 \beqlb
 \label{generatorY}
{\mbb A}_mf(y_1,\cdots,y_m)&=& G^mf(y_1,\cdots,y_m)\cr
&&+\frac{1}{2}\sum_{1\leq i\neq j\leq
m}\left(f(\theta_{ij}(y_1,\cdots,y_m))-f(y_1,\cdots,y_m)\right),
 \eeqlb
where $f\in C^2(\dfR^m)$ and $\theta_{ij}(x_1,\cdots,x_m)$ denote
the element of $\dfR^m$ obtained by replacing $x_j$ by $x_i$ in
$(x_1,\cdots,x_m)$.
\end{lemma}
 {\bf Proof. } The proof is similar to that of part (b) in Lemma 2.1 of
 \cite{[DK96]}.  For $1\leq i,\,j\leq m$,
 define
 \begin{align}\label{tildeB}
 \tlb_{j0}&=B_{\alpha 0},&\text{where
 }\alpha&=\Sigma_j^m(0),&\cr
 \tlb_{ijk}&= B_{\alpha\beta\gamma},
 &\text{where
 }\alpha&=\Sigma_i^m(\tilde{\tau}_{ijk}-),
 \beta=\Sigma_j^m(\tilde{\tau}_{ijk}-),& \cr
 &&\qquad\gamma&=N_{\alpha \beta}(\tilde{\tau}_{ijk}-)&
 \end{align}
 Define
$\tilde{\gamma}_{ijk}=\min\{\tilde{\tau}_{i'jk'},\, i'\neq j:
\tilde{\tau}_{i'jk'}>\tilde{\tau}_{ijk}\}$ and let
$\tilde{\gamma}_{j0}$ be the first jump time of
$\tilde{N}_j\equiv\sum_{i\neq j}\tilde{N}_{ij}$. By Lemma
\ref{Aselection}, $Y_j^m(t)=X_{\Sigma^m_j(t)}^m(t)$ yields that for
$0\leq
 t<\tilde{\gamma}_{j0}$
 \beqlb
 \label{Y1}
 Y^m_j(t)=Y^m_j(0)
 +\ez \tlb_{j0}(t)+\int_0^t\int_{\dfR}h(y-Y^m_{j}(s))W(dyds)
 \eeqlb
and for $\tilde{\tau}_{ijk}\leq
 t<\tilde{\gamma}_{ijk}$,
 \beqlb
 \label{Y2}
 Y^m_j(t)=Y_{i}(\tilde{\tau}_{ijk})+\ez (\tlb_{ijk}(t)
 -\tlb_{ijk}({\tilde{\tau}_{ijk}}))
 +\int_{\tilde{\tau}_{ijk}}^t\int_{\dfR}h(y-Y^m_{j}(s))W(dyds).
 \eeqlb
 By Lemmas A5.1 and A5.2 of \cite{[DK96]},  $\{\tlb_{j0}\}, \{\tlb_{ijk}\}$
 and $\{Y_j(0)\}$ are independent of $\{\tilde{N}_{ij}\}$ and
 $\Sigma^m$. Furthermore,
 the $\tlb_{j0}$ and the $\tlb_{ijk}$ are independent Brownian motions
  and $(Y_1^m(0),\cdots,Y^m_m(0))$ has the same distribution
 as $(X_1(0),\cdots,X_m(0))$. Then the desired result follows
 from (\ref{Y1}) and (\ref{Y2}). \qed

By (\ref{Y1}) and (\ref{Y2}), we see $(Y_1^m(t),\cdots, Y^m_m(t))$
is exchangeable and has the same empirical measures as
$(X_1,\cdots,X_m).$ From the construction above, $\Sigma^m(t)$ must
be independent of $Y^m(t)$. Thus for each $t>0$,
$(X_1(t),X_2(t),\cdots)$ is exchangeable. To show the existence of
$(\L, \mu)$-martingale problem, we need the following lemma.
 \begin{lemma}\label{Lemma2.3}
\begin{enumerate}
 \item[(a).] Suppose that  $Z(t)$ is a  $P(\dfR)$-valued process
 satisfying the martingale formula (\ref{martfor}) for every $F\in \D(\L)$.
  Then $\{Z(t):t\geq0\}$ has a continuous modification and
  for $\phi\in C^2(\dfR)$
 \beqlb\label{l2.3a}
 M_t(\phi):=\la Z(t),\phi\ra-\la Z(0),\phi\ra-\frac{\rho_{\ez}}{2}\int_0^t\la Z(s),\phi'' \ra ds
 \eeqlb
is a martingale with quadratic variation
 \beqlb\label{l2.3b}
 \gamma\int_0^t \left(\la Z(s),\phi^2\ra-\la Z(s),\phi\ra^2\right)ds+
 \int_0^t ds\int_{\dfR}\la Z(s),h(\cdot-y)\phi'\ra^2dy.
 \eeqlb
\item[(b).] If a continuous $P(\dfR)$-valued process $Z(t)$
satisfies the martingale problem (\ref{l2.3a}) and (\ref{l2.3b}),
then it is also a solution of $(\L, \mu)$-martingale problem.
\end{enumerate}
 \end{lemma}
{\bf Proof. } (a). The existence of continuous modification follows
from Lemma 2.1 of \cite{[EK87]} and the fact that (\ref{martfor}) is
a martingale for each $F\in \D(\L)$ which also yields (\ref{l2.3a})
and (\ref{l2.3b}). The proof for assertion (b) is a classical
approximation procedure. We left it to the interested readers. \qed

Now, we come to our main result in this section.

\begin{theorem}
\label{ThExist}  Given $\mu\in P(\mbb R)$, suppose that
$\{X_i(0),i\geq1\}$ is an exchangeable sequence of random variables
such that
$$\lim_{m\rar\infty}\frac{1}{m}\sum_{i=1}^m
 \dz_{X_i(0)}=\mu.$$ Let
 \beqlb\label{approdirac}
 Z_m(t)=\frac{1}{m}\sum_{i=1}^m
 \dz_{X_i(t)}=\frac{1}{m}\sum_{i=1}^m\dz_{Y_i^m(t)}.\eeqlb Then
 the $(\L,\mu)$-martingale problem has a solution $Z$ such that for each
 $t>0$,
 \beqlb
 \label{ThExista}
 \lim_{m\rar\infty}\sup_{s\leq t}\rho(Z_m(s),Z(s))=0\quad a.s.,
 \eeqlb
where $\rho$ denotes the Prohorov metric on $P(\mbb R)$.
\end{theorem}
{\bf Proof. } With the help of Lemma \ref{Lemma2.3} which can be
regarded as a version of Lemma 2.3 of \cite{[DK96]}, the proof is
similar to Theorem 2.4 of \cite{[DK96]}. We omit it here. \qed

\section{Sample Path Properties}\label{SECstate}
In this section, we show that when $\ez>0$, $Z(t)$ is absolutely
continuous respect to $dx$ for almost all $t\geq0$ and when $\ez=0$
the values of $Z$ are purely atomic. We first describe the
\textsl{weak atomic topology} on $M(\mbb R)$ introduced by Ethier
and Kurtz \cite{[EK94]}. Recall that $\rho$ denotes the Prohorov
metric on $M(\mbb R)$, which induces the topology of the weak
convergence. Define the metric $\rho_a$ on $M(\mbb R)$ by
 \beqlb
 \label{atomictopo}
 \rho_a(\mu,\nu)=\rho(\mu,\nu)&+&\sup_{0<\ez\leq1}\bigg{|}\int_{\dfR}
 \int_{\dfR}\Phi(|x-y|/\ez)\mu(dx)\mu(dy)\cr
 \ar\ar\qquad\quad-\int_{\dfR}
 \int_{\dfR}\Phi(|x-y|/\ez)\nu(dx)\nu(dy)\bigg{|},
 \eeqlb
where $\Phi(\cdot)=\left(1-\cdot\right)_+.$ The topology on $M(\mbb
R)$ induced by $\rho_a$ is called the \textsl{weak atomic topology}.
For $\mu\in M(\mbb R)$, define $\mu^*=\sum\mu(\{x\})^2\dz_{x}$. We
need the following results of \cite{[EK94]}.
\begin{lemma}
\label{EKlemmas} Let $\mu_n,\mu\in M(\dfR)$. \begin{enumerate}
\item[(a).] Suppose $\rho(\mu_n,\mu)\rar0$. Then
$\rho(\mu_n^*,\mu^*)\rar0$ if and only if
$\mu_n^*(\dfR)\rar\mu^*(\dfR)$;
\item[(b).] $\rho_a(\mu_n,\mu)\rar0$ if and only if
$\rho(\mu_n,\mu)\rar0$ and $\rho(\mu_n^*,\mu^*)\rar0$;
\item[(c).] Suppose $Z\in C_{(M(\dfR),\,\rho)}[0,\infty)$.
If $Z^*(\dfR)\in C_{[0,\infty)}[0,\infty)$, then $Z\in
C_{(M(\dfR),\,\rho_a)}[0,\infty)$.
\end{enumerate}
\end{lemma}
{\bf Proof.} See Lemmas 2.1, 2.2 and 2.11 of \cite{[EK94]} for (a),
(b)  and (c), respectively.\qed

 Our
first main result in this section is the following theorem.

\begin{theorem}
\label{Thstate} Suppose $Z$ is a solution of $(\L,\mu)$-martingale
problem. Assume $\ez=0$.   Then ${\bf P}\{Z(t)\in P_{a}(\mbb
R),\,t>0\}={\bf P}\{Z(\cdot)\in
C_{(M(\dfR),\,\rho_a)}[0,\infty)\}=1,$ where $P_a(\mbb R)$ denotes
the collection of purely atomic probability measures on $\mbb R$.
\end{theorem}
{\bf Proof.}  According to the look down construction,
(\ref{lookdown1}) and (\ref{lookdown2}), if $X_j$ `looks down'
$X_i$, and assume the value of $X_i$ at the jump time, then  $X_j$
and $X_i$ have the same sample path before the next jump time.
Define
 $$
x_i(t)=X_i(0)+\int_0^t\int_{\mbb R}h(y-x_i(s))W(dyds),\quad
t\geq0,\quad i=1,2,\cdots.
 $$
Therefore, by Lemma \ref{EKlemmas}, $Z_m(\cdot)\in D_{(P(\mbb R),\,
\rho_a)}[0,\infty) $ and $Z^*_m(t,\mbb R)$ is monotone in $t\geq0$.
According to Proposition 3.3 of \cite{[DK96]} and Lemma
\ref{lehomo}, almost surely for $t>0$, there are only finite number
paths, denoted by $D(t)$ which is independent of $m$, alive in the
`look down system'. Let ${t_0}>0$ be fixed. Note that  $D$ is
c\`{a}dl\`{a}g on $ [t_0,+\infty)$. Typically, $D(t)\leq D(s)$ for
$t>s$. Let $\{x_{c_i}(t_0), i=1,2,\cdots, D(t_0) \}$ be the
enumeration of the living paths at $t_0$ with
$x_{c_1}(t_0)<x_{c_2}(t_0)<\cdots<x_{c_{D(t_0)}}(t_0)$. Thus for
$t>t_0$, we may represent $Z_m(t)$ by
 \beqlb\label{Thstatec}
Z_m(t)=\sum_{i=1}^{D(t_0)}\frac{b_{i,m}(t)}{m}\dz_{x_{c_i}(t)},\quad
t\geq t_0,
 \eeqlb
 where $b_{i,m}(t), i=1,2,\cdots$ are nonnegative integer-valued
c\`{a}dl\`{a}g random processes defined on $ [t_0,+\infty)$ with
$\sum_{i=1}^{D(t)} b_{i,m}=m$.   Note that by Lemma \ref{lehomo},
for every $T>t_0$, almost surely, \beqlb\label{Thstated}\inf_{i\neq
j}\inf_{t_0\leq t\leq T}|x_{c_i}(t)-x_{c_j}(t)|>0.\eeqlb Therefore,
according to (\ref{ThExista})  we may represent $Z(t)$ by
 \beqlb\label{repZ}
Z(t)=\sum_{i=1}^{D(t_0)}{b_{i}(t)}\dz_{x_{c_i}(t)},\quad t\geq t_0,
 \eeqlb
 where $b_{i}(t)\geq 0, i=1,2,\cdots$ are
c\`{a}dl\`{a}g random processes defined on $ [t_0,+\infty)$ with
 \beqlb\label{Thstatee}
\sup_{t_0\leq t\leq
T}\sum_{i=1}^{D(t_0)}|b_{i,m}(t)/m-b_{i}(t)|\rar0,\quad a.s.\quad
\text{as }m\rar\infty.
 \eeqlb
 Since $t_0$ is arbitrary, ${\bf
P}\{Z(t)\in P_{a}(\mbb R),\,t>0\}=1.$  From above and Lemma
\ref{EKlemmas}, we see $Z(\cdot\vee t_0)\in D_{(P(\mbb
 R),\,\rho_a)}[0,\infty),\,a.s.$ Typically,
$$
Z_m^*(\cdot\vee t_0,\mbb R)\rar Z^*(\cdot\vee t_0,\mbb R)\quad
\text{in}\quad D_{\mbb
 R}[0,\infty)\quad \text{as}\quad m\rar\infty\quad a.s.
$$
On the other hand, according to the `look down construction',
 if we define
$$
J(Z_m^*(t\vee t_0,\dfR)):=\int_0^{\infty}e^{-u}[1\wedge\sup_{0\leq
t\leq u}|Z_m^*(t\vee t_0,\dfR)-Z_m^*((t\vee t_0)-,\dfR)|]du,
$$
then
$$
J(Z_m^*(t\vee t_0,\dfR))\leq \frac{4m+2}{m^2}\rar 0\quad\text{as
}m\rar\infty.
$$
By Theorem 3.10.2 of \cite{[EK86]} and Lemma \ref{EKlemmas},
$Z(\cdot\vee t_0)\in C_{(P(\mbb
 R),\,\rho_a)}[0,\infty),\,a.s.$
 Set  $D=\{(x,y)\in\dfR^2:x=y\}$ and $D_2=D\times \dfR^2+\dfR^2\times
D$. By approximating an indicate function from continuous functions,
we see that ($\ref{Dual}$) holds for $f=1_{D}$ and $g=1_{D_2}$. Note
that $\la Z(t)^2, f\ra=Z^*(t,\dfR)$ and $\la Z(t)^2, f\ra^2=\la
Z(t)^4,g\ra$. Therefore, by (\ref{Dual}) and the right continuity of
$(F_t,M_t)$,
$$
\lim_{t\downarrow0}{\bf
E}|Z^*(t,\dfR)-\mu^*(\dfR)|^2=\lim_{t\downarrow0}{\bf E}|\la Z(t)^2,
f\ra-\la \mu^2,f\ra|^2=0.
$$
By Lemma \ref{EKlemmas} and the monotonicity of $Z^*_m(t,\mbb R)$,
$\rho_a(Z(t),\mu)\rar0$ almost surely as $t\rar0$. Thus $Z(\cdot)\in
C_{(P(\mbb
 R),\,\rho_a)}[0,\infty),\,a.s.$\qed

In the next theorem, we shall show that when $\ez>0$ $Z(t,dx)$ is
absolutely continuous with respect to $dx$ and derive the SPDE for
the density.
\begin{theorem}\label{Thcon}
Suppose $Z$ is a solution of $(\L,\mu)$-martingale problem. Assume
$\ez>0$. Then for $t>0$, $Z(t,dx)$ is absolutely continuous with
respect to $dx$ and the density $Z_t(x)$ satisfies the following
SPDE: for $\phi\in {\cal S}(\mbb R)$,
 \beqlb\label{SPDE}
 \la
 Z_t,\phi\ra-\la\mu,\phi\ra\ar=\ar\int_0^t\int_{\dfR}\sqrt{\gamma Z_s(x)}\phi(x){V}(dsdx)
 -\int_0^t\int_{\dfR}\la Z_s,\phi\ra\sqrt{\gamma Z_s(x)}{V}(dsdx)\cr
\ar\ar+\int_0^t\int_{\mbb R}\la Z_s,h(x-\cdot)\phi'\ra
W(dsdx)+\frac{\rho_{\ez}}{2}\int_0^t\la Z_s,\phi''\ra ds,
 \eeqlb
 where $V$ and $W$ are two independent Brownian sheets and ${\cal S}(\dfR)$
 is the space of rapidly decreasing $C^{\infty}$-function defined on
 $\dfR$ equipped with the Schwartz topology.
\end{theorem}
{\bf Proof.} We borrow the ideas in Theorem 1.7 of \cite{[KS88]}.
First by dual relationship (\ref{Dual}), one can derive that for any
$\phi,\psi\in C(\mbb R)$,
 \beqlb\label{moment1}
 \E\la Z(t),\phi\ra=\la\mu, T_t^1\phi\ra
 \eeqlb
and
 \beqlb\label{moment2}
 \E\left[\la Z(t),\phi\ra\la Z(t),\psi\ra\right]=e^{-\gamma t}
\la \mu^2,T_t^2\phi\psi\ra+ \int_0^t e^{-\gamma s}\la \mu
T_{t-s}^1,\Psi_{12}(T_s^2
 \phi\psi)\ra ds.
 \eeqlb
For $\epsilon>0$, the semigroup $(T_t^m)_{t>0}$ is uniformly
elliptic and has density $q_m(t,x,y)$ satisfying
$$q_m(t,x,y)\leq c\cdot g_m({\ez'
t},x,y),~~t>0,~x,y\in\mathbb{R}^m,$$ where $c$ is a constant and
$g_m(t,x,y)$ denotes the transition density of the $m$-dimensional
standard Brownian motion; see \cite{[Dy65]}. Without loss of
generality, we assume $\ez'=1$.
  Note that
 \beqnn
 &&\int_{\dfR^2}q_1(u,x,z_1)q_1(u',x,z_2)q_1(t-s,z,y)
 q_2(s,(y,y),(z_1,z_2))dz_1dz_2\\
&\rar& q_1(t-s,z,y)
 q_2(s,(y,y),(x,x))
 \eeqnn
as $u,u'\rar0$. Meanwhile, \beqnn
 &&\int_{\dfR^2}q_1(u,x,z_1)q_1(u',x,z_2)q_1(t-s,z,y)
 q_2(t,(y,y),(z_1,z_2))dz_1dz_2\\
&&\leq c\int_{\dfR^2}g_1(u,x,z_1)g_1(u',x,z_2)g_1(t-s,z,y)
 g_2(t,(y,y),(z_1,z_2))dz_1dz_2\\
 &&=cg_1(u+s,x,y)g_1(u'+s,x,y)g_1(t-s,z,y).
 \eeqnn
Take $\phi=\phi_{u,x}=q_1(u,x,\cdot)$ and
$\psi=\psi_{u',x}=q_1(u',x.\cdot)$ in (\ref{moment2}). By dominated
convergence theorem, when $u,u'\rar0$,
 \beqlb\label{Thcona}
 &&\int_0^Tdt\int dx \int_0^t e^{-\gamma s}\la \mu
T_{t-s}^1,\Psi_{12}(T_s^2
 \phi\psi)\ra ds \cr
 &&\rar\int_0^Tdt\int dx\int_0^tds\int_{\dfR^2}e^{-\gamma s}
q_1(t-s,z,y)
 q_2(s,(y,y),(x,x))dy\mu(dz).
 \eeqlb
Similarly, we have
 \beqlb\label{Thconb}
 &&\int_0^Tdt\int dx  e^{-\gamma t}
\la \mu^2,T_t^2\phi\psi\ra\cr
 &&\rar\int_0^Tdt\int dx\int_{\dfR^4}e^{-\gamma t}
 q_2(t,(x_1,x_2),(x,x))\mu(dx_1)\mu(dx_2).
 \eeqlb
Combining (\ref{Thcona}) and (\ref{Thconb}) together yields
  $\{\la Z(t), q_{u}(x,\cdot)\ra, u>0\}$ is a Cauchy
  sequence in $L^2(\Omega\times[0,T]\times \dfR)$.
This implies the existence of the density $Z_t(x)$ of $Z_t$ in
$L^2(\Omega\times[0,T]\times \dfR).$

\smallskip

Next, we derive the SPDE (\ref{SPDE}). Choose an one dimensional
standard Brownian motion $\hat{B}_t$ independent of $Z_t$. For any
fixed $c>1/2$, set $G_t=\exp(\hat{B}_t+(c-1/2)t)$. So $Z_t>0$ and
$Z_t\rar\infty$ as $t\rar\infty$ a.s. It also satisfies
$$
dG_t=\sqrt{\gamma}G_td\hat{B}_t+cG_tdt,\quad G_0=0.
$$
Define $C_t=\int_0^t G_sds$. $C_t$ is strictly increasing and
$C_t\rar\infty$ as $t\rar\infty$ a.s.. Let $C_t^{-1}$ denote its
inverse function on $[0,\infty).$ Define measure-valued process
$I_t$ by
$$
I_t(dx)=G_{C_t^{-1}}\cdot Z_{C_t^{-1}}(dx).
$$
By Ito's formula, (\ref{l2.3a}) and (\ref{l2.3b})
 \beqnn
\la I_t,\phi\ra&=&\la I_0,\phi\ra+\int_0^{C_t^{-1}}G_sdM_s(\phi)
 +\int_0^{C_t^{-1}}\sqrt{\gamma} G_s\la Z_s,\phi\ra d\hat{B}_s\\
 &&~+c\int_0^{C_t^{-1}}G_s\la Z_s,\phi\ra ds+
 \frac{\rho_{\ez}}{2}\int_0^{C_t^{-1}}G_s\la Z_s,\phi''\ra ds.
 \eeqnn
Then $$\tilde{M}_t(\phi):=\int_0^{C_t^{-1}}G_sdM_s(\phi)
 +\int_0^{C_t^{-1}} G_s\la Z_s,\phi\ra d\hat{B}_s,
\quad t\geq0,
 $$is a local martingale with quadratic function
$$
\la \tilde{M}(\phi)\ra_t={\gamma}\int_0^t\la I_s,\phi^2\ra
ds+\int_0^tds \int_{\mbb R}\la I_s,h(x-\cdot)\phi'\ra^2/\la I_s,1\ra
dx.
$$
Clearly, $I_t(dx)$ is also absolutely continuous with respect to
$dx$. Denote the corresponding density by $I_t(x)$. Similar to the
martingale representation theorem (see Theorem 3.3.6 of
\cite{[KX95]} or Theorem III-7 of \cite{[EM90]}),  there exists two
independent $L^2(\mbb R)$-cylindrical Brownian motion $\tilde{V}$
and $\tilde{W}$ (may be on an extension probability space) such that
$$
\tilde{M}_t(\phi)=\int_0^t\la
f(s,I_s)^*\phi,d\tilde{V}_s\ra_{L^2(\mbb R)}+\int_0^t\la
g(s,I_s)^*\phi,d\tilde{W}_s\ra_{L^2(\mbb R)},
$$
where $f(s,I_s)$ and $g(s,I_s)$ are linear maps from $L^2(\mbb R)$
to $\S'(\mbb R)$, the space of Schwartz distributions,  such that
for $\phi\in \S(\dfR)$,
$$f(s,I_s)^*\phi(x)=\sqrt{{\gamma}I_s(x)}\phi(x)$$ and
$$
g(s,I_s)^*\phi(x)=\int_{\mbb R}h(x-y)\phi'(y)I_s(y)dy/\sqrt{\la
I_s,1\ra}.
$$ Thus
\beqnn \la I_t,\phi\ra&=&\int_0^t\la
f(s,I_s)^*\phi,d\tilde{V}_s\ra_{L^2(\mbb R)}+\int_0^t\la
g(s,I_s)^*\phi,d\tilde{W}_s\ra_{L^2(\mbb R)}\\
 &&~+c\int_0^t\la I_s,\phi\ra/\la I_s,1\ra ds+
 \frac{\rho_{\ez}}{2}\int_0^t\la I_s,\phi''\ra/\la I_s,1\ra ds.
 \eeqnn
Define two new $L^2(\mbb R)$-cylindrical Brownian motions $\hat{V}$
and $\hat{W}$ by
$$
\la \hat{V}_t,\phi\ra=\int_0^{C_t}\frac{1}{\la I_s,1\ra}\la
d\tilde{V},\phi\ra,\quad\la
\hat{W}_t,\phi\ra=\int_0^{C_t}\frac{1}{\la I_s,1\ra}\la
d\tilde{W},\phi\ra.
$$
Since $\tilde{V}$ and $\tilde{W}$ are independent, $\hat{V}$ and
$\hat{W}$ are orthogonal (hence they are independent).
 Then we can find two independent
Brownian sheets ${V}(dtdx)$ and ${W}(dtdx)$  such that
$$
\tilde{V}_t(l)=\int_0^t\int_{\dfR}l(x){V}(dsdx),\quad
\tilde{W}_t(l)=\int_0^t\int_{\dfR}l(x){W}(dsdx),\quad \forall\, l\in
L^2(\dfR).
$$
Using Ito's formula and noting that $\la Z_t,\phi\ra=\la
I_{C_t},\phi\ra/\la I_{C_t},1\ra$ yield \beqnn
 \la
 Z_t,\phi\ra-\la\mu,\phi\ra\ar=\ar\int_0^t\int_{\dfR}\sqrt{\gamma Z_s(x)}\phi(x){V}(dsdx)
 -\int_0^t\int_{\dfR}\la Z_s,\phi\ra\sqrt{\gamma Z_s(x)}{V}(dsdx)\cr
\ar\ar+\int_0^t\int_{\mbb R}\la Z_s,h(x-\cdot)\phi'\ra
W(dsdx)+\frac{\rho_{\ez}}{2}\int_0^t\la Z_s,\phi''\ra ds
 \eeqnn
for $\phi\in {\cal S}(\mbb R)$. We have completed the proof. \qed

\section{Connections to Measure-valued Branching Processes in a Random Medium}
It has been shown that there are deep connections between the
Dawson-Watanabe and Fleming-Viot superprocesses; see
\cite{[EM91],{[KS88]},[P92]}. In this section, we shall show that
the Fleming-Viot processes in random environment is a class of
measure-valued branching processes in a Brownian medium, conditioned
to have total mass one. Such measure-valued branching processes were
first constructed and studied by \cite{[W97]} and \cite{[W98]}. The
argument in this section is similar to those in \cite{[P92]} with
some modifications. Let $\{\omega(t),t\geq0\}$ and
$\{\hat{\omega}(t),t\geq0\}$ denote the coordinate processes on
$C_{P(\dfR)}[0,\infty)$ and $C_{M(\dfR)}[0,\infty)$, respectively.
Define ${\cal F}^0_t=\sigma(\omega(s); s\leq t)$, $\hat{\cal
F}^0_t=\sigma(\hat{\omega}(s); s\leq t)$, ${\cal F}_t={\cal
F}_{t+}^0$ and $\hat{\cal F}_t=\hat{\cal F}_{t+}^0$. Based on the
results in \cite{[W98]} and the continuity of $\hat{\omega}$, for
each $\mu\in M(\dfR)$, there exists an unique probability measure
$\hat{\bf Q}_{\mu}$ on $C_{M(\dfR)}[0,\infty)$ such that for
$\phi\in C^2(\dfR)$
 \beqlb\label{SDSM1}
 \hat{M}_t(\phi):=\la \hat{\omega}
 (t),\phi\ra-\la\mu,\phi\ra-\frac{\rho_{\ez}}{2}\int_0^t\la \hat{\omega}(s),\phi'' \ra
 ds, \quad t\geq0,
 \eeqlb
under $\hat{\bf Q}_{\mu}$ is a continuous $\hat{\cal
F}_t$-martingale starting at 0 with quadratic variation
 \beqlb\label{SDSM2}
 \la\hat{M}(\phi)\ra_t=\gamma\int_0^t \la \hat{\omega}(s),\phi^2\ra ds+
 \int_0^t ds\int_{\dfR}\la \hat{\omega}(s),h(\cdot-y)\phi'\ra^2dy.
 \eeqlb
 Let
 \beqnn
 C_+&=&\{f:[0,\infty)\rar[0,\infty): f \textrm{ continuous }, \exists\,t_f\in(0,\infty]
\textrm{ such that }\\
&&\qquad f(t)>0 \textrm{ if }t\in[0,t_f) \textrm{ and }f(t)=0
\textrm{ if } t\geq t_f\}
 \eeqnn
with the compact-open topology.
 Let $L_y\in P(C_+)$ denote
the law of the unique solution of
 $$
 \eta_t=y+\int_0^t\sqrt{\gamma\eta_s}dB_s,
 $$
where $B$ is a standard Brownian motion. Note that
 \beqlb\label{Thdwfv2}
 \hat{\bf Q}_{\mu}(\hat{\omega}(\dfR)\in\cdot)=L_{\mu(\dfR)}(\cdot).
 \eeqlb
For $\mu\in M(\dfR)-\{0\}$, define
$\bar{\mu}(\cdot)=\mu(\cdot)/\mu(\dfR)$. Let $\{{\bf
Q}_{\bar{\mu},f}(A):A\in {\cal F},f\in C_+\}$ be a regular
conditional probability for $\bar{\omega}$ given
$\hat{\omega}_{\cdot}=f(\cdot)$ under $\hat{\bf Q}_{\mu}$, where
$\cal F $ denotes the Borel $\sigma$-field on
 $C_{P(\dfR)}[0,\infty)$. That is
 $$
 \hat{\bf Q}_{\mu}(\bar{\omega}\in
 A|\hat{\omega}_{\cdot}(\dfR)=f(\cdot))={\bf Q}_{\bar{\mu},f}(A)
 \quad \forall A\in {\cal F}.
 $$\begin{lemma}\label{ThDWFV}
  For each $\mu\in M(\dfR)-\{0\}$, there exists a subset $C_{\mu}$ of $C_+$ such that
  $L_{\mu(\dfR)}(C_{\mu})=1$ and for $f\in C_{\mu}$,
 under ${\bf Q}_{\bar{\mu},f}$
\beqlb\label{ThDWFVm0}
 {M}_t^f(\phi,\omega):=\la {\omega}_t,
 \phi\ra-\la\bar{\mu},\phi\ra-\frac{\rho_{\ez}}{2}\int_0^t\la
{\omega}_s,\phi'' \ra
 ds, \quad t<t_f,
 \eeqlb
is an ${\cal F}_t$-martingale starting at 0 for every $\phi\in
C^2(\dfR)$ with
 \beqlb
 \label{ThDWFVm1}
 \la{M}^f(\phi)\ra_t&=&\gamma\int_0^{t}(\la{\omega}_s,\phi^2\ra-
 \la{\omega}_s,\phi\ra^2)f(s)^{-1}ds
 \cr&&\quad+\int_0^{t} ds\int_{\dfR}\la
 {\omega}_s,h(\cdot-y)\phi'\ra^2dy\quad \forall\,t<t_f \eeqlb
 and $\omega_t=\omega_{t_f} \textrm{ for all }t\geq t_f$.
 \end{lemma}
 \begin{remark}
 Note that if $f=1$, then (\ref{ThDWFVm0}) and (\ref{ThDWFVm1}) are
 just (\ref{l2.3a}) and (\ref{l2.3b}), respectively.
 \end{remark}
{\bf Proof.} Define $T_n=\inf\{t:\hat{\omega}_t(\dfR)\leq 1/n\}$
 and for $\phi\in C^2(\dfR)$\beqlb
 \label{ThDWFVa}
 \bar{M}_t^n(\phi):=\int_0^{t\wedge T_n}
 \hat{\omega}_s(\dfR)^{-1}d\hat{M}_s(\phi)-\int_0^{t\wedge T_n}
 \la\hat{\omega}_s,\phi\ra\hat{\omega}_s(\dfR)^{-2}d\hat{M}_s(1).
 \eeqlb
Thus for fixed $t$, $\{\bar{M}_t^n(\phi):n\geq1\}$ is a martingale
in $n$. By Ito's formula,
 \beqlb
 \label{ThDWFVb}
 \la\bar{\omega}_{t\wedge T_n},\phi\ra
 =\la\bar{\mu},\phi\ra+\frac{\rho_{\ez}}{2}\int_0^{t\wedge
 T_n}\la\bar{\omega}_s,\phi''\ra ds+\bar{M}_t^n(\phi),
 \eeqlb
which implies that
 \beqlb
 \label{ThDWFVc}
 \sup_{t\leq K,\, n\geq 1}|\bar{M}_t^n(\phi)|\leq
 2||\phi||_{\infty}+\frac{K\rho_{\ez}}{2}||\phi''||_{\infty}.
 \eeqlb
Therefore, according to the Martingale Convergence Theorem and
maximal inequality,  $\bar{M}_t^n(\phi)$ converges as $n\rar \infty$
uniformly for $t$ in compacts a.s. (by perhaps passing to a
subsequence). We denote by $\bar{M}_t(\phi)$ the limit  which is a
continuous martingale satisfying
 \beqlb
 \label{ThDWFVd}
 \bar{M}_t^n(\phi)=\bar{M}_{t\wedge T_n}(\phi),\quad \forall
 t\geq0,\quad a.s.
 \eeqlb
and \beqlb
 \label{ThDWFVe}
 \sup_{t\leq K}|\bar{M}_t(\phi)|\leq
 2||\phi||_{\infty}+\frac{K\rho_{\ez}}{2}||\phi''||_{\infty}.
 \eeqlb
Letting $n\rar\infty$ in (\ref{ThDWFVb}) yields \beqlb
 \label{ThDWFVf}
 \la\bar{\omega}_{t},\phi\ra
 =\la\bar{\mu},\phi\ra+\frac{\rho_{\ez}}{2}\int_0^{t\wedge
 T_0}\la\bar{\omega}_s,\phi''\ra ds+\bar{M}_t(\phi),\quad
 \forall t\geq0~a.s.~\forall \phi\in C^2(\dfR),
 \eeqlb
where $T_0=\inf\{t:\hat{\omega}_t(\dfR)=0\}$. Note that
 \beqlb
 \label{ThDWFV1}
 \bar{M}_{t\wedge T_0}(\phi)=\bar{M}_{t}(\phi).
 \eeqlb
 Let $s<t$
and let $F$ be a bounded
$\sigma(\hat{\omega}_{\cdot}(\dfR))$-measurable random variable.
Since $\{\hat{\omega}_t(\dfR):t\geq0\}$ is a martingale under
$\hat{\bf Q}_{\mu}$, the martingale representation theorem implies
that there exists some $\sigma(\hat{\omega}_s(\dfR):s\leq t)$-
predictable function $f$ such that
 \beqlb\label{ThDWFVg}
 F=\hat{\bf
 Q}_{\mu}(F)+\int_0^{\infty}f(s,\hat{\omega})d\hat{\omega}_s(\dfR).
 \eeqlb
According to (\ref{ThDWFVd}) and (\ref{ThDWFVg}),
 \beqnn
\ar\ar \hat{\bf Q}_{\mu}((\bar{M}_{t\wedge
T_n}(\phi)-\bar{M}_{s\wedge T_n}(\phi))F|{\cal F}_s)\\
\ar\ar\quad=\hat{\bf
Q}_{\mu}((\bar{M}^n_{t}(\phi)-\bar{M}^n_{s}(\phi))
\int_0^{\infty}f(s,\hat{\omega})d\hat{\omega}_s(\dfR)|{\cal F}_s)\\
\ar\ar\quad =\hat{\bf Q}_{\mu}((\int_{s\wedge T_n}^{t\wedge T_n}
 \hat{\omega}_u(\dfR)^{-1}d\hat{M}_u(\phi)-\int_{s\wedge T_n}^{t\wedge T_n}
 \la\hat{\omega}_u,\phi\ra\hat{\omega}_u(\dfR)^{-2}d\hat{M}_u(1))
\int_s^tf(s,\hat{\omega})d\hat{\omega}_u(\dfR)|{\cal F}_s)\\
\ar\ar\quad=\hat{\bf Q}_{\mu}(\int_{s\wedge T_n}^{t\wedge
T_n}(\la\hat{\omega}_u,\phi\ra\hat{\omega}_u(\dfR)^{-1}
-\la\hat{\omega}_u,\phi\ra\hat{\omega}_u(\dfR)^{-1})f(u)du|{\cal
F}_s)\\
\ar\ar\quad=0
 \eeqnn
By letting $n\rar\infty$ in the above, we have
 $$
\hat{\bf Q}_{\mu}((\bar{M}_t(\phi)-\bar{M}_s(\phi))F|{\cal F}_s)=0,
 $$
which yields for a fixed $\phi\in C^2(\dfR)$,
$\{\bar{M}_t(\phi):t\geq0\}$ is a martingale with respect to ${\cal
G}_t:={\cal F}_t\vee\sigma(\hat{\omega}_s(\dfR):s\geq0)$. On the
other hand, by (\ref{ThDWFVd}) and (\ref{ThDWFV1}),
 \beqlb
 \label{ThDWFVh}
 \la\bar{M}(\phi)\ra_t&=&\gamma\int_0^{t\wedge T_0}(\la\bar{\omega}_s,\phi^2\ra-
 \la\bar{\omega}_s,\phi\ra^2)\hat{\omega}_s(\dfR)^{-1}ds
 \cr&&\quad+\int_0^{t\wedge T_0} ds\int_{\dfR}\la
 \bar{\omega}_s,h(\cdot-y)\phi'\ra^2dy\qquad\hat{\bf Q}_{\mu}-a.s.
 \eeqlb
Set $M_t^f(\phi,\omega)=M_{t_f-}^f(\phi)$ for $t\geq t_f$. By
(\ref{ThDWFVf}) and (\ref{ThDWFV1}),
 \beqlb\label{ThDWFVi}
 \bar{M}_t(\phi)=M_t^{\omega_{\cdot}(\dfR)}(\phi,\bar{\omega}),\quad
 \forall\, t\geq0\quad\hat{\bf Q}_{\mu}-a.s.~\forall \phi\in C^2(\dfR).
 \eeqlb
Then for each $G\in b{\cal F}_t^0$ and $s<t$, by the ${\cal G}_t$
martingale property of $\bar{M}_t(\phi)$, (\ref{Thdwfv2}) and
(\ref{ThDWFVi}),
$$
{\bf
Q}_{\bar{\mu},f}\left(\left(M_t^f(\phi)-M_s^f(\phi)\right)G\right)=0\quad
L_{\mu(\dfR)}-a.a.f.
$$
By considering rational and the fact that $C_{P(\dfR)}[0,\infty)$
with local uniform topology is a standard measurable space and
taking limits in $s$ and $G$, we could find a $L_{\mu(\dfR)}$-null
set off which the above holds for all $s<t$ and $G\in {\cal F}_s$.
That is $\{M_t^f(\phi):t\geq0\}$ is an ${\cal F}_t$-martingale under
${\bf Q}_{\bar{\mu},f}$ for $L_{\mu(\dfR)}-a.a.f.$ Take
$t_n^f=\inf\{u:f(u)\leq 1/n\}$. According to (\ref{ThDWFVh}) and
above arguments, we can deduce that for every $n\geq1$
$$
M_{t\wedge t_n^f}^f(\phi)^2-\gamma\int_0^{t\wedge
t_n^f}(\la{\omega}_s,\phi^2\ra-
 \la{\omega}_s,\phi\ra^2)f(s)^{-1}ds
 -\int_0^{t\wedge
t_n^f} ds\int_{\dfR}\la
 {\omega}_s,h(\cdot-y)\phi'\ra^2dy,\quad t\geq0,
$$
is an ${\cal F}_t$-martingale under ${\bf Q}_{\bar{\mu},f}$ for
$L_{\mu(\dfR)}-a.a.f.$ Now, consider a countable subset of
$C^2(\dfR)$, $C_S(\dfR)$, such that we can approximate any function
$\phi\in C^2(\dfR)$ by a sequence $\{\phi_k:\,k\geq1\}\subset
C_S(\dfR)$ in such a way that not only $\phi$ but all of its
derivatives up to the second order are approximated boundedly and
pointwise. Taking limits in $M_t^f(\phi)$ and $\la M^f(\phi)\ra_t$
yields the desired conclusion.   \qed

\bigskip

For $T>0$, define $(\Omega_{T-},{\cal F}_{T-})=(C_{P({\dfR})}[0,T),
\text{Borel sets}).$ $(\hat{\Omega}_{T-},\hat{\cal F}_{T-})$ denotes
the same space with $M(\dfR)$ in place of $P(\dfR)$.  If ${\bf Q}$
is a probability on $C_{P(\dfR)}[0,\infty),$ then ${\bf Q}|_{T-}$ is
defined on $(\Omega_{T-},{\cal F}_{T-})$ by ${\bf Q}|_{T-}(A)={\bf
Q}(\omega|_{[0,T)}\in A).$ Similarly, one defines $(\Omega_{T},{\cal
F}_{T})$, $(\hat{\Omega}_{T},\hat{\cal F}_{T})$ and ${\bf Q}|_{T}.$
Suppose ${\bf Q}_{\mu}$ is the unique probability measure on
$C_{P(\dfR)}[0,\infty)$ such that $\{\omega(t),t\geq0\}$ under $\bf
Q_{\mu}$ is a solution of $(\cal L,\mu)$-martingale problem. Our
main result in this subsection is the following theorem which is
analogous to Corollary 4 of \cite{[P92]}.

\begin{theorem}\label{ThM}
Suppose that $\{\mu_n\}\subset M(\dfR)-\{0\}$ satisfy
$\bar{\mu}_n\rar\mu$ in $P(\dfR)$.
\begin{enumerate}
\item[(a).] If for each $n$, there exists a function $f_n\in
C_{\mu_n}$ such that for some $T>0$, $\sup_{0\leq t\leq
S}|f_n-1|\rar 0$ for $S<T$ as $n\rar\infty$, then
\beqlb\label{Thm1}{\bf Q}_{\bar{\mu}_n,f_n} |_{T-}\rar {\bf
Q}_{\mu}|_{T-}\text{ weakly on } (\Omega_{T-},{\cal F}_{T-}).\eeqlb
\item[(b).] Let $\{A_n\}$ be a sequence of Borel subset of $C_+$ such that
$ L_{\mu_n(\dfR)}(A_n)>0$ for every $n\geq1$. If for some $T>0$
$$
\sup\{|g(t)-1|:g\in A_n,t\leq S\}\rar0 \text{ as
}n\rar\infty,\forall S<T,
$$
then $$ \hat{\bf Q}_{\mu_n}(\bar{\omega}\in
\cdot|\omega_{\cdot}(\dfR)\in A_n)|_{T-}\rar {\bf
Q}_{\mu}|_{T-}\text{ weakly on } (\Omega_{T-},{\cal F}_{T-}). $$
\end{enumerate}
\end{theorem}
{\bf Proof.} (a). It suffices to prove
$$
{\bf Q}_{\bar{\mu}_n,f_n} |_{S}\rar {\bf Q}_{\mu}|_{S}\text{ weakly
on } (\Omega_{S},{\cal F}_{S}).
$$
Let $\hat{\mathbb{R}}=\dfR\cup\{\partial\}$ denote the  one-point
compactification of $\mathbb{R}$. Since $\sup_{0\leq t\leq
S}|f_n-1|\rar 0$ for $S<T$ as $n\rar\infty$, $\inf_{t\leq S}f_n\geq
1/2$ for $n$ larger enough and
 $$
 |\la M^{f_n}(\phi)\ra_t-\la M^{f_n}(\phi)\ra_s|\leq
 \frac{\gamma}{2}||\phi||_{\infty}^2|t-s|
 +||\rho||_{\infty}||\phi'||^2_{\infty}|t-s|,\quad \forall\, s,t\leq S,
  {\bf Q}_{\bar{\mu}_n,f_n}-a.s.
 $$
By Theorem 2.3 of \cite{[RC86]}, one can check that $\{{\bf
Q}_{\bar{\mu}_n,f_n}|_S: n\geq1\}$ is tight in
$P(C_{P(\hat{\dfR})}[0,S])$. Let ${\bf Q}$ be a limit point in
$P(C_{P(\hat{\dfR})}[0,S])$. With abuse of notation, we denote by
$\{\omega_s:s\leq S\}$  the coordinate processes of
$C_{P(\hat{\dfR})}[0,S]$. One may use Skorohod representation
theorem to see that under $\bf Q$,
 \beqlb\label{l2.3a1}
 M_t(\phi):=\la \omega_t,\phi\ra-\la \mu,\phi\ra-\frac{\rho_{\ez}}
 {2}\int_0^t\la \omega_s,\phi'' \ra ds
 \eeqlb
is a continuous martingale starting at 0 for $t\leq S$ and $\phi\in
C^2_{\partial}({\mathbb{R}})$ with quadratic variation
 \beqlb\label{l2.3b2}
 \gamma\int_0^t \left(\la \omega_s,\phi^2\ra-\la \omega_s,\phi\ra^2\right)ds+
 \int_0^t ds\int_{\dfR}\la \omega_s,h(\cdot-y)\phi'\ra^2dy.
 \eeqlb
  We claim that
$$\mathbf{Q}\{\omega_t(\{\partial\})=0~\textrm{for all}~
 t\in[0,S]\}=1.$$ Consequently, $\mathbf{Q}$ is supported by
 $C_{P({\dfR})}[0,S]$.
For $k\geq1$, let
$$
 \phi_k(x)=\begin{cases}
 \exp\{-\frac{1}{|x|^2-k^2}\},& \textrm{if}~|x|>k,\\
  0,&\textrm{if}~|x|\leq k.
  \end{cases}
  $$
One can check that $\{\phi_k\}\subset C^2_{\partial}({\mathbb{R}})$
such that $\lim_{|x|\rightarrow\infty}\phi_k(x)=1$,
$\lim_{|x|\rightarrow\infty}\phi_k(x)'=0$ and
$\phi_k(\cdot)\rightarrow1_{\{\partial\}}(\cdot)$ boundedly and
pointwise. $||\phi_k'||\rightarrow0$ and $||\phi_k''||\rightarrow0$
as $k\rightarrow\infty$.   By martingale inequality, we have
\begin{eqnarray*}
&&\mathbf{Q}\{\sup_{0\leq t\leq S}|M_t(\phi_k)-M_t(\phi_j)|^2\}\cr
&&\quad\leq4\gamma\int_0^S\mathbf{Q}_{\mu}\{\la\omega_s,(\phi_k-\phi_j)^2\ra\}
ds+8\gamma\int_0^S\mathbf{Q}_{\mu}\{\la\omega_s,|\phi_k-\phi_j|\ra\}
ds\cr&&\qquad+4\int_0^Sds\int_{\hat{\mathbb{R}}}\mathbf{Q}\{\la
 \omega_s,h(z-\cdot)(\phi_k'-\phi_j')\ra^2\}dz.\end{eqnarray*}
 By dominated convergence theorem, $\mathbf{Q}\{\sup_{0\leq
 t\leq
S}|M_t(\phi_k)-M_t(\phi_j)\ra|^2\}\rightarrow0$ as
$k,j\rightarrow\infty$. Therefore, there exists
$M^{\partial}=(M^{\partial}_t)_{t\leq S}$ such that for every $t\leq
S$,
$$\mathbf{Q}\{|M_t(\phi_k)-M_t^{\partial}|^2\}
\rightarrow0
$$
and (by perhaps passing to a subsequence)
$$
\sup_{0\leq s\leq t}|M_s(\phi_k)-M_s^{\partial}| \rightarrow0\quad
{\bf Q}-a.s.
$$
as $k\rightarrow\infty$.  We obtain $M^{\partial}$ is a continuous
martingale. It follows from ($\ref{l2.3a1}$) that
$M_t^{\partial}=\omega_t(\{\partial\})$ is a continuous martingale
with mean zero . Thus $\mathbf{Q}(\omega_t(\{\partial\}))=0$. Then
the claim follows from the  continuity of
$\big{\{}\omega_t(\{\partial\}):t\geq0\big{\}}$. Extend $\bf Q$ to
$C_{P(\dfR)}[0,\infty)$ by setting the conditional distribution of
$\{\omega_{t+S}:t\geq0\}$ given ${\cal F}_S^0$ equal to ${\bf
Q}_{\omega_S}$. Then ${\bf Q}={\bf Q}_{\mu}$ and so ${\bf Q}|_S={\bf
Q}_{\mu}|_S$. We complete the proof of (a).

 \bigskip

 (b). Let $H:\Omega|_{T-}\rar \dfR$ be bounded and continuous. Then
 by Lemma \ref{ThDWFV} and (\ref{Thm1}),
  \beqnn
  &&|\hat{\bf Q}_{\mu_n}(H(\bar{\omega})
  |\omega_{\cdot}(\dfR)\in A_n)- {\bf Q}_{\mu}(H)|\cr
 &&\quad = |\hat{\bf Q}_{\mu_n}(H(\bar{\omega})
  |\omega_{\cdot}(\dfR)\in A_n\cap C_{\mu_n})- {\bf Q}_{\mu}(H)|\\
  && \quad\leq \left|\int_{A_n\cap C_{\mu_n}}
{\bf Q}_{\bar{\mu}_n,g}(H)-{\bf
Q}_{\mu}(H)dL_{\mu_n(\dfR)}L_{\mu_n(\dfR)}(A_n\cap C_{\mu_n})^{-1}\right|\\
&&\quad \leq \sup_{g\in A_n\cap C_{\mu_n}}|{\bf
Q}_{\bar{\mu}_n,g}(H)-{\bf Q}_{\mu}(H)|\\
&&\quad \rar0\quad \text{ as }\quad n\rar\infty.
  \eeqnn
We are done. \qed
\begin{corollary}
\label{coroPE} Suppose that $\{\mu_n\}\subset M(\dfR)-\{0\}$ satisfy
$\bar{\mu}_n\rar\mu$ in $P(\dfR)$. For $T>0$, let $T_n\rar T$ and
$\dz_n\rar0$ and assume $|\mu_n(\dfR)-1|<\dz_n$. Then
\begin{enumerate}
\item[(a).] $\hat{\bf Q}_{\mu_n}(\bar{\omega}\in \cdot|
\sup\limits_{t\leq T_n}|\omega_t(\dfR)-1|<\dz_n)
\xrightarrow{weakly} {\bf Q}_{\mu}|_{T-}\quad \text{on}\quad
(\Omega_{T-},{\cal F}_{T-})$;
\item[(b).] $\hat{\bf Q}_{\mu_n}(\hat{\omega}\in \cdot|
\sup\limits_{t\leq T_n}|\omega_t(\dfR)-1|<\dz_n)
\xrightarrow{weakly} {\bf Q}_{\mu}|_{T-}\quad \text{on}\quad
(\hat{\Omega}_{T-},\hat{\cal F}_{T-})$.
\end{enumerate}
\end{corollary}
{\bf Proof.} Setting
$$
A_n=\{g\in C_+:\sup_{t\leq T_n}|g(t)-1|<\dz_n\}
$$
and Theorem \ref{ThM} yield (a). (b) follows from (a) and the fact
that for $S<T$ and $n$ large enough, $$ \hat{\bf
Q}_{\mu_n}\left(\sup_{t\leq
S}|\hat{\omega}_t(\dfR)^{-1}-1|<\frac{\dz_n}{1-\dz_n}\bigg{|}A_n\right)=1.
$$
\qed

\bigskip

\textbf{Acknowledgement}. I would like to give my sincere thanks to
Professors Shui Feng, Zenghu Li, Jie Xiong and Hao Wang for their
simulating discussions.

\bigskip

\centerline{APPENDIX}
\appendix
\section{Random selections of stochastic integrals}
\begin{lemma}
\label{Aselection} Let $W(dsdy)$ be a space-time white noise on
$[0,\infty)\times \dfR$ based on Lebesgue measure measure. Let
$\{X_i(t),t\geq0,i=1,2, \cdots\}$ be a sequence of real valued
predictable stochastic processes. Let $h(x,y)$ be a measurable
function on $\dfR\times\dfR$. Define stochastic integrals
$$
Y_i(t):=\int_0^t\int_{\dfR}h(X_i(s),y)W(dsdy),\quad t\geq0,\quad
i\geq 1.
$$
Suppose $\pi$ is a random variable taking values in
$\{1,2,\cdots\}$, independent of $\{X_i,i=1,2,\cdots\}$ and $W$.
Then
$$
Y_{\pi}(t)=\int_0^t\int_{\dfR}h(X_{\pi}(s),y)W(dsdy),\quad t\geq0.
$$
\end{lemma}
{\bf Proof.} If $h$ is a simple function, the desired conclusion is
obvious. For general result, one can consider the $L^2$
approximation and Ito's isometry; see Theorem 2.2.5 of
\cite{[Wl86]}. \qed

\section{Stochastic flow of diffeomorphism}
In this part, we consider the following stochastic differential
equation \beqlb\label{SDEQN}
\xi_t=x+\int_s^t\int_{\dfR}h(y-\xi(s))W(dsdy),\quad x\in\dfR,\quad
t\geq s,
 \eeqlb
where $W(dsdy)$ is a space-time white noise on $[0,\infty)\times
\dfR$ based on Lebesgue measure measure. The existence and pathwise
uniqueness for (\ref{SDEQN}) have been proved in \cite{[DLW01]}.

\begin{lemma}
\label{lehomo} Suppose $h\in C^2(\mbb R)$. There is a modification
of the solution, denoted by $\xi_{s,t}(x)$, such that almost surely
\begin{enumerate}
\item[(1)]  $\xi_{s,t}(x,\omega)$ is
continuous in $(s,t,x)$ and satisfies $\lim_{t\downarrow
s}\xi_{s,t}(x,\omega)=x;$

\item[(2)]
$\xi_{s,t+u}(x,\omega)=\xi_{t,t+u}(\xi_{s,t}(x,\omega),\omega)$ is
satisfied for all $s<t$ and $u>0$;

\item[(3)]   the map
$\xi_{s,t}(\cdot,\omega):\dfR\rar\dfR$ is an onto homeomorphism for
all $s<t$;

\item[(4)]  the map
$\xi_{s,t}(\cdot,\omega):\dfR\rar\dfR$ is a $C^2$-diffeomorphism for
all $s<t$.
\end{enumerate}

\end{lemma}
{\bf Proof. } The argument is exactly similar to that in Chapter 2
of \cite{[Ku84]}. We omit it here and left it to interested readers.
\qed

\end{document}